\documentclass[12pt, reqno]{amsart}
\usepackage{latexsym, amsfonts, amsmath, amssymb}

\newtheorem{theorem}{Theorem}[section]

\newtheorem{lemma}[theorem]{Lemma}
\newtheorem{proposition}[theorem]{Proposition}
\newtheorem{corollary}[theorem]{Corollary}

\newtheorem*{MT}{Main Theorem}

\theoremstyle{definition}
\newtheorem{remark}[theorem]{Remark}

\newtheorem{example}[theorem]{Example}

\newcommand{\abs}[1]{\ensuremath{\left\vert #1 \right\vert}}
\newcommand{\Cal}[1]{\ensuremath{\mathcal{#1}}}

\def\C{{\mathbb C}}

\def\N{{\mathbb N}}
\def\Z{{\mathbb Z}}
\def\Q{{\mathbb Q}}

\def\O_K{{\Cal{O}_{K}}}
\def\O_F{{\Cal{O}_{F}}}
\def\N_F{{\Cal{N}_{F/\Q}}}

\def\F{{\mathbb{F}}}
\def\deg{{\textrm{deg}}}
\def\ord{{\textrm{ord}}}
\def\ker{{\textrm{ker}}}

\numberwithin{equation}{section}
\numberwithin{theorem}{section}

\begin{document}

\title{multiple zeta values over global function fields}

\author{Riad Masri}

%\address{Department of Mathematics, University of Texas, Austin, TX USA 78712}

\address{Max-Planck-Institut f\"ur Mathematik, Vivatsgasse 7, 53111 Bonn, Germany}

\email{masrirm@mpim-bonn.mpg.de}

\subjclass{Primary: 11M41; Secondary: 14H05}
\keywords{Function field; Multiple zeta function; Meromorphic continuation; Rational function.}

\begin{abstract}
Let $K$ be a global function field with finite constant field $\mathbb{F}_{q}$ of order $q$. 
In this paper we develop the analytic theory of a multiple zeta function $Z_{d}(K; s_{1}, \ldots, s_{d})$ in $d$ 
independent complex variables defined over $K$. This is the function field analog of the Euler-Zagier 
multiple zeta function $\zeta_{d}(s_{1},\ldots, s_{d})$ of depth $d$ (\cite{Z1}). 
Our main result is that  $Z_{d}(K;s_{1}, \ldots, s_{d})$ has a meromorphic 
continuation to all $(s_{1}, \ldots, s_{d})$ in $\C^{d}$ and is a rational function in 
each of $q^{-s_{1}}, \ldots, q^{-s_{d}}$ with a specified denominator.
\end{abstract}

\maketitle

\section{Introduction and statements of results}

1. In \cite{Z1} Zagier defined the multiple zeta function of depth $d$ by
\begin{align*}
\zeta_{d}(s_{1}, \ldots, s_{d})=\sum_{0 < n_{1} < \cdots < n_{d}}n_{1}^{-s_{1}}\cdots n_{d}^{-s_{d}},
\quad (s_{1}, \ldots, s_{d}) \in \C^{d},
\end{align*}
which is absolutely convergent and analytic in the region
\begin{align}\label{region}
\mathrm{Re}(s_{k}+\cdots +s_{d})>d-k+1, \quad k=1,\ldots, d.
\end{align}
He then defined the multiple zeta values of depth $d$ by  
\begin{align*}
\zeta_{d}(a_{1},\ldots, a_{d})=\sum_{0 < n_{1} < \cdots < n_{d}}n_{1}^{-a_{1}}\cdots n_{d}^{-a_{d}}, \quad 
a_{k} \in \Z_{\geq 1}, \quad a_{d}>1,
\end{align*}
and described in part the fundamental role these numbers play in geometry, number theory, and physics. 
This has continued to be revealed in the work of Drinfeld, Goncharov, Kontsevich, 
Manin, and Zagier, among others (see the discussion below).

In this paper we begin to develop a theory of multiple zeta values over global function fields. 
Our objective here is to develop the analytic properties of two new multiple zeta functions analogous to 
$\zeta_{d}$. The first of these is associated to the polynomial 
ring $\F_{q}[T]$ over a finite field $\F_{q}$ of order $q$, and the second of these 
is associated to a global function field $K$ with finite constant field $\F_{q}$. 
We will prove that each of these functions has a meromorphic continuation to all 
$(s_{1}, \ldots, s_{d})$ in $\C^{d}$ and is a rational function in each of $q^{-s_{1}}, \ldots, q^{-s_{d}}$
with a specified denominator.  

We now describe the main result of this paper. Let $K$ be a global function field with finite constant field
$\mathbb{F}_{q}$. For a divisor $D$ of $K$, let $\textrm{deg}(D)$ be its degree and
$\abs{D}=q^{\textrm{deg}(D)}$ be its norm. Let 
$\Cal{D}_{K}^{+}$ be the sub semi-group of effective divisors of $K$. We 
define the multiple zeta function of depth $d$ over $K$ by
\begin{align*}
Z_{d}\left(K;s_{1}, \ldots, s_{d}\right)=
\sum_{\substack{(D_{1}, \ldots, D_{d}) \in \Cal{D}_{K}^{+} \times \cdots \times \Cal{D}_{K}^{+}\\
0 \leq \deg(D_{1}) \leq \cdots \leq \deg(D_{d})}}\prod_{k=1}^{d}\abs{D_{k}}^{-s_{k}},
\end{align*} 
which is absolutely convergent and analytic in the region (\ref{region}).
Our main result is the following 

\begin{MT}
The multiple zeta function $Z_{d}\left(K;s_{1},\ldots, s_{d}\right)$  
has a meromorphic continuation to all $(s_{1},\ldots,s_{d})$ in $\C^{d}$ and is a rational function in 
each of $q^{-s_{1}}, \ldots, q^{-s_{d}}$ with a specified denominator. Further, $Z_{d}\left(K;s_{1},\ldots, s_{d}\right)$  
has possible simple poles on the linear subvarieties 
\begin{align*}
s_{k}+\cdots +s_{d}=0,1, \ldots, d-k+1, \quad k=1, \ldots, d.
\end{align*}
\end{MT}

The proof of the main theorem rests on the Riemann-Roch theorem for global function fields. 

In the remaining part of the introduction we want to provide some motivating background on the multiple
zeta values over $\Q$ and describe the results of this paper in more detail.
As mentioned, the multiple zeta values and their generalizations, the multiple polylogarithms at $N$-th roots of unity,
\begin{align*}
\mathrm{Li}_{a_{1}, \ldots, a_{d}}\left(\zeta_{N}^{\alpha_{1}}, \ldots, \zeta_{N}^{\alpha_{d}}\right)=
\sum_{0 < n_{1} < \cdots < n_{d}}\frac{\left(\zeta_{N}^{\alpha_{1}}\right)^{n_{1}}
\cdots \left(\zeta_{N}^{\alpha_{d}}\right)^{n_{d}}}
{n_{1}^{a_{1}}\cdots n_{d}^{a_{d}}}, \quad \zeta_{N}=e^{2 \pi i /N},
\end{align*}
have been the focus of much attention in the past 15 years. In an important development,   
Kontsevich discovered that the multiple zeta values can be expressed as an iterated integral of  
Chen-type. This led to A. Goncharov's 
interpretation of the multiple zeta values as periods of mixed motives, and his remarkable work 
\cite{G1, G2, G3, G4, G5} on mixed Tate motives over $\textrm{Spec}(\Z)$ and proof of the 
upper bound in Zagier's dimension conjecture. 

Zagier's dimension conjecture is perhaps the central problem in the subject. Recall that this 
is a statement about the non-trivial $\Q$-linear relations which arise
between multiple zeta values of the same weight (see \cite{Z1}). 
Let $Z_{w}$ denote the $\Q$-algebra generated by all multiple zeta values of weight 
$w=a_{1}+\cdots +a_{d}$.
The dimension conjecture states that $\dim_{\Q}Z_{w}=d_{w}$, where $d_{0}=1,~
d_{1}=0,~d_{2}=1,$ and $d_{w}=d_{w-2}+d_{w-3}$ for $w \geq 3$. For $w\geq 2$, this
formula implies that $d_{w}$ is less than or equal to $2^{w-2}$, the number of multiple 
zeta values of weight $w$, and hence that there are 
many non-trivial $\Q$-linear relations between multiple zeta values of the same weight.  
Goncharov proved the inequality $\dim_{\Q}Z_{w} \leq d_{w}$. 

Additional interesting work on the multiple zeta values and related topics can be found in  
Drinfeld \cite{D}, Goncharov and Manin \cite{GM}, Kontsevich \cite{K1, K2}, Manin \cite{Ma1, Ma2}, 
and Zagier \cite{Z2}. The particularly nice survey article of Kontsevich and Zagier \cite{KZ} discusses
the multiple zeta values in the more general context of periods and special values of $L$-functions.

There has also been a great deal of attention focused on the \textit{analytic} properties of 
multiple zeta functions. For example, the meromorphic continuation of $\zeta_{d}$ to all $(s_{1},\ldots, s_{d})$ in $\C^{d}$ was
first established by Goncharov and Kontsevich in \cite{G4}. The existence of such a meromorphic 
continuation was originally obscured by the presence of points of indeterminacy, which are special 
types of singularities arising in several complex variables (see \cite{M}). Goncharov and Kontsevich established 
the meromorphic continuation by using a $d$-dimensional Mellin transform to express  
$\zeta_{d}$ as the pairing of a meromorphic distribution and a test function in a 
certain modified Schwartz class. These methods were subsequently extended by J. Kelliher and the 
author in \cite{KM} to give a sufficient condition for the meromorphic continuation of a more 
general class of multiple Dirichlet series of Euler-Zagier type. 
\vspace{0.05in}

2. Let $\F_{q}[T]$ be the ring of 
polynomials with coefficients in the finite field $\F_{q}$ of order $q$. For $f$ in $\F_{q}[T]$, 
let $\textrm{deg}(f)$ be its degree and $\abs{f}=q^{\textrm{deg}(f)}$.
Recall that the zeta function of $\F_{q}[T]$ is defined by
\begin{align*}
Z(\F_{q}[T],s)=
\sum_{\substack{f \in \F_{q}[T]\\ f~\textrm{monic}}}\abs{f}^{-s}, \quad \textrm{Re}(s)>1. 
\end{align*}
We formally define the multiple zeta function of depth $d$ over $\F_{q}[T]$ by
\begin{align}\label{MDSFiniteFields}
Z_{d}\left(\F_{q}[T];s_{1},\ldots, s_{d}\right)=
\sum_{\substack{(f_{1},\ldots, f_{d}) \in \F_{q}[T] \times \cdots \times \F_{q}[T]\\ 
f_{1},\ldots, f_{d}~\textrm{monic}\\
0 \leq \textrm{deg}(f_{1}) \leq \cdots \leq \textrm{deg}(f_{d})}}
\prod_{k=1}^{d}\abs{f_{k}}^{-s_{k}}.
\end{align}

The following result establishes the analytic properties of (\ref{MDSFiniteFields}).

\begin{theorem}\label{MainTheorem1}

(1) $Z_{d}\left(\F_{q}[T];s_{1},\ldots, s_{d}\right)$ is absolutely 
convergent and analytic in the region (\ref{region}).

(2) $Z_{d}\left(\F_{q}[T];s_{1},\ldots, s_{d}\right)$ has a 
meromorphic continuation to all $(s_{1},\ldots, s_{d})$ in $\C^{d}$ and is a rational function 
in each of $q^{-s_{1}}, \ldots, q^{-s_{d}}$. In fact,
\begin{align*}
\prod_{k=1}^{d}\left(1-q^{d-k+1-(s_{k}+\cdots +s_{d})}\right)
Z_{d}\left(\F_{q}[T];s_{1},\ldots, s_{d}\right)=1.
\end{align*}

(3) $Z_{d}\left(\F_{q}[T];s_{1},\ldots, s_{d}\right)$ has simple poles on the 
linear subvarieties
\begin{align*}
s_{k}+\cdots +s_{d}=d-k+1, \quad k=1,\ldots, d.
\end{align*}

(4) The function 
\begin{align*}
\xi_{d}\left(\F_{q}[T];s_{1},\ldots, s_{d}\right):=
\prod_{k=1}^{d}\frac{q^{-(s_{k}+\cdots +s_{d}-(d-k))}}{1-q^{-(s_{k}+\cdots +s_{d}-(d-k))}}
Z_{d}\left(\F_{q}[T];s_{1},\ldots, s_{d}\right)
\end{align*}
is invariant under the involution
\begin{align*}
s_{1} \mapsto 2d-1-2(s_{2}+\cdots +s_{d})-s_{1}.
\end{align*}

(5) $Z_{d}\left(\F_{q}[T];s_{1}, \ldots, s_{d}\right)$ has the Euler product
\begin{align*}
Z_{d}\left(\F_{q}[T];s_{1}, \ldots, s_{d}\right)=
\prod_{\substack{P \in \F_{q}[T]\\ P~\mathrm{monic}\\P~\mathrm{irreducible}}}
\prod_{k=1}^{d}\left(1-\frac{1}{\abs{P}^{s_{k}+\cdots +s_{d}-(d-k)}}\right)^{-1}.
\end{align*}
\end{theorem}

\begin{remark} 
$Z_{d}\left(\F_{q}[T];s_{1}, \ldots, s_{d}\right)$ seems to be the first example of 
a Zagier type multiple zeta function possessing a functional equation and Euler product. In fact, it is 
not difficult to see that $\zeta_{d}$ cannot satisfy a functional equation in any of the independent
variables $s_{1}, \ldots, s_{d}$. This is because of a shifting in the variables $n_{1}, \ldots, n_{d}$ which 
occurs if one expresses $\zeta_{d}$ as a sum over $\Z^{d}_{\geq 1}$.
\end{remark}

\begin{remark}
It is possible to compute the residues at the simple poles on the linear subvarieties
\begin{align*}
s_{k}+\cdots +s_{d}=d-k+1, \quad k=1, \ldots, d.
\end{align*}
For example, suppose $d=2$. A straightforward application of L'Hopital's rule shows that the 
double zeta function  $Z_{2}\left(\F_{q}[T];s,w\right)$
has simple poles on the linear subvarieties $s+w=2$ and $w=1$ with the following residues:
\begin{enumerate}
\item For $s\neq 1$, 
\begin{align*}
\lim_{w \rightarrow 1}(w-1)Z_{2}\left(\F_{q}[T];s,w\right)=
\frac{\left(1-q^{1-s}\right)^{-1}}{\log(q)}.
\end{align*}

\item For $w \neq 1$, 
\begin{align*}
\lim_{s \rightarrow 2-w}(s-(2-w))Z_{2}\left(\F_{q}[T];s,w\right)=
\frac{\left(1-q^{1-w}\right)^{-1}}{\log(q)}.
\end{align*}

\item For $s \neq 1$, 
\begin{align*}
\lim_{w \rightarrow 2-s}(w-(2-s))Z_{2}\left(\F_{q}[T];s,w\right)=
\frac{\left(1-q^{s-1}\right)^{-1}}{\log(q)}.
\end{align*}
\end{enumerate}
\end{remark}

From Theorem \ref{MainTheorem1} we obtain the following 
factorization of $Z_{d}\left(\F_{q}[T];s_{1}, \ldots, s_{d}\right)$.

\begin{corollary}\label{FactorizationCorollary}
$Z_{d}\left(\F_{q}[T];s_{1}, \ldots, s_{d}\right)$ has the factorization
\begin{align*}
Z_{d}\left(\F_{q}[T];s_{1}, \ldots, s_{d}\right)=
\prod_{k=1}^{d}Z\left(\F_{q}[T],s_{k}+\cdots +s_{d}-(d-k)\right).
\end{align*}
\end{corollary}

\begin{proof}
This follows from Theorem \ref{MainTheorem1}, part $(2)$, and the identity
$Z(\F_{q}[T],s)=1/(1-q^{1-s})$.
\end{proof}

\begin{corollary}
$Z_{d}\left(\F_{q}[T];s_{1}, \ldots, s_{d}\right)$ has no zeros in $\C^{d}$.
\end{corollary}

\begin{proof}
This follows from Corollary \ref{FactorizationCorollary} and the fact that $Z(\F_{q}[T],s)$ has
no zeros in $\C$.
\end{proof}

\begin{remark}
Upon examining the proof of part (4) of Theorem \ref{MainTheorem1}, one sees that
the functional equation is actually satisfied by
\begin{align*}
\frac{q^{-(s_{1}+\cdots +s_{d}-(d-1))}}{1-q^{-(s_{1}+\cdots +s_{d}-(d-1))}}
Z_{d}\left(\F_{q}[T];s_{1},\ldots, s_{d}\right).
\end{align*} 
However, as suggested by Corollary \ref{FactorizationCorollary}, by including the factors
\begin{align*}
\frac{q^{-(s_{k}+\cdots +s_{d}-(d-k))}}{1-q^{-(s_{k}+\cdots +s_{d}-(d-k))}}, \quad k=2, \ldots, d, 
\end{align*}
we can obtain additional functional relations. These are no longer involutions, but 
instead involve a mixing of the variables. For example, suppose $d=2$. Let $$w \mapsto 1-w$$ in 
$$\xi_{2}\left(\F_{q}[T];s,w\right).$$
Then arguing as in the proof of part (4) of Theorem \ref{MainTheorem1}, we obtain the functional relation
\begin{align*}
\xi_{2}\left(\F_{q}[T];s,1-w\right)
&=\xi_{2}\left(\F_{q}[T];s-2w+1,w\right).
\end{align*}
\end{remark}
\vspace{0.05in}

3. We now recall some background on function fields.
A function field in one variable over a constant field $F$ is a field $K$ 
containing $F$ and at least one element $x$ transcendental over $F$ such that $K/F(x)$ 
is a finite algebraic extension. A function field in one variable over a finite constant field
$F=\F_{q}$ is called a global function field. Throughout this paper we assume that $K$ is a global
function field. 

A prime in $K$ is a discrete valuation ring $R$ with maximal ideal $P$ such that $\F_{q} \subset P$ and
the quotient field of $R$ equals $K$. The degree $\deg(P)$ of $P$ is the dimension of $R/P$ over $\F_{q}$,
which is finite. The group $\Cal{D}_{K}$ of divisors of $K$ is 
the free abelian group generated by the primes in $K$. A typical divisor $D$ is written additively by 
$D=\sum_{P}a(P)P.$ The degree of $D$ is defined by $\deg(D)=\sum_{P}a(P)\deg(P).$

Given $a \in K^{*}$, the divisor $(a)$ of $a$ is defined by
$(a)=\sum_{P}\ord_{P}(a)P.$ The map $K^{*} \rightarrow \Cal{D}_{K}$ defined by $a \mapsto (a)$
is a homomorphism whose image $\Cal{P}_{K}$ is the group of principal divisors.
Two divisors $D_{1}$ and $D_{2}$ are linearly equivalent $D_{1}\sim D_{2}$ if their difference is 
principal; that is, $D_{1}-D_{2}=(a)$ for some $a \in K^{*}$. Define the divisor class group by
$Cl_{K}=\Cal{D}_{K}/\Cal{P}_{K}$. 

It can be shown that the degree of a principal divisor is zero 
(see \cite{R}, Proposition 5.1). Thus, the degree map $\deg:Cl_{K} \rightarrow \Z$ is a homomorphism.
Let $\ker(\deg)=Cl_{K}^{0}$ be the group of divisor classes of degree zero. 
It can also be shown that the number $\abs{Cl_{K}^{0}}$ of divisor classes of degree zero is finite 
(see \cite{R}, Lemma 5.6). Define the class number of $K$ to be $h_{K}=\abs{Cl_{K}^{0}}$.
Because $K$ has divisors of degree one (see \cite{S}), we obtain the exact sequence
\begin{align*}
0 \rightarrow Cl_{K}^{0} \rightarrow Cl_{K} \rightarrow \Z \rightarrow 0.
\end{align*}

A divisor $D$ is effective if $a(P) \geq 0$ for all $P$. We denote this by $D \geq 0$.
Given a divisor $D$, let 
\begin{align*}
L(D)=\left\{x \in K^{*}:~(x)+D \geq 0\right\} \cup \{0\}.
\end{align*}
It can be shown that $L(D)$ is a finite dimensional vector space over $\F_{q}$. Let $l(D)$
the dimension of $L(D)$ over $\F_{q}$. 

We are now in a position to state the following form of the Riemann-Roch theorem for global function fields. 

\begin{theorem}[Riemann-Roch]\label{RiemannRoch}
There is an integer $g \geq 0$ and a divisor class $\Cal{C}$ such that for $C \in \Cal{C}$
and $A \in \Cal{D}_{K}$ we have
\begin{align*}
l(A)=\mathrm{deg}(A)-g+1+l(C-A).
\end{align*}
\end{theorem}

The integer $g$, which is uniquely determined by $K$, is called the genus of $K$.
\vspace{0.05in}

4. For a divisor $D$ of $K$, let $\abs{D}=q^{\deg(D)}$ be its norm. 
Then $\abs{D}$ is a positive integer,
and for any two divisors $D_{1}$ and $D_{2}$, $\abs{D_{1}+D_{2}}=\abs{D_{1}}\abs{D_{2}}$.
Let $\Cal{D}_{K}^{+}$ be the sub semi-group of effective divisors of $K$. 
Recall that the zeta function of $K$ is defined by 
\begin{align*}
Z(K,s)=\sum_{D \in \Cal{D}_{K}^{+}}\abs{D}^{-s}, \quad \textrm{Re}(s)>1.
\end{align*}
We formally define the multiple zeta function of depth $d$ over $K$ by
\begin{align}\label{MDSFunctionFields}
Z_{d}\left(K;s_{1}, \ldots, s_{d}\right)=
\sum_{\substack{(D_{1}, \ldots, D_{d}) \in \Cal{D}_{K}^{+} \times \cdots \times \Cal{D}_{K}^{+}\\
0 \leq \deg(D_{1}) \leq \cdots \leq \deg(D_{d})}}\prod_{k=1}^{d}\abs{D_{k}}^{-s_{k}}.
\end{align}

If $K$ is a global function field of genus $g=0$, then $K=\F_{q}(T)$ is the rational function field. 
The following result establishes the analytic properties of (\ref{MDSFunctionFields}) in the case $K=\F_{q}(T)$.

\begin{theorem}\label{MainTheorem2} 
(1) $Z_{d}\left(\F_{q}(T);s_{1},\ldots, s_{d}\right)$ is absolutely 
convergent and analytic in the region (\ref{region}).

(2) $Z_{d}\left(\F_{q}(T);s_{1},\ldots, s_{d}\right)$ has a meromorphic continuation 
to all $(s_{1},\ldots,s_{d})$ in $\C^{d}$ and is a rational function in each of 
$q^{-s_{1}}, \ldots, q^{-s_{d}}$.
In fact,
\begin{align*}
Q\left(q^{-s_{1}},\ldots, q^{-s_{d}}\right)
Z_{d}\left(\F_{q}(T);s_{1},\ldots, s_{d}\right)
\end{align*}
is a polynomial of degree $ \leq 2d-1$ in each of $q^{-s_{1}}, \ldots, q^{-s_{d}}$, 
where
\begin{align*}
&Q\left(q^{-s_{1}},\ldots, q^{-s_{d}}\right)=
(q-1)^{d}\left(1-q^{-s_{d}}\right)\left(1-q^{1-s_{d}}\right)\\
& \times \prod_{k=1}^{d-1}\left(1-q^{-(s_{k}+\cdots +s_{d})}\right)
\left(1-q^{1-(s_{k}+\cdots +s_{d})}\right)
\left(1-q^{2-(s_{k}+\cdots +s_{d})}\right).
\end{align*}

(3) $Z_{d}\left(\F_{q}(T);s_{1},\ldots, s_{d}\right)$ has possible simple poles on 
the linear subvarieties 
\begin{align*}
s_{k}+\cdots +s_{d}=0,1, \ldots, d-k+1, \quad k=1, \ldots, d.
\end{align*}
\end{theorem}

Recall that a multiple zeta value is said to be reducible (completely reducible) if it can be written as a rational
linear combination of products of lower depth (depth one) multiple zeta values.
The following corollary shows that $Z_{d}\left(\F_{q}(T);s_{1}, \ldots, s_{d}\right)$ 
is always a rational linear combination of products of 1-dimensional zeta functions over 
$\F_{q}[T]$.

\begin{corollary}\label{DecompositionCorollary} 
$Z_{d}\left(\F_{q}(T);s_{1}, \ldots, s_{d}\right)$ is a 
rational linear combination of products of zeta functions of the type
\begin{align*}
Z\left(\F_{q}[T],s_{k}+\cdots +s_{d}-1\right),~Z\left(\F_{q}[T],s_{k}+\cdots +s_{d}\right),~
Z\left(\F_{q}[T],s_{k}+\cdots +s_{d}+1\right),
\end{align*}
where $k=1, \ldots, d$.
\end{corollary}

Corollary \ref{DecompositionCorollary} is in stark contrast to what occurs for the multiple zeta values,
where one must place restrictions on the depth $d$ and weight $w$ to guarantee reducibility.
For example, when the depth $d$ and weight $w$ of a multiple zeta value have different parity, 
the multiple zeta value is reducible (see \cite{T}). 
As a first instance of this, there is the following fact due to Euler and Zagier:
\textit{Every double zeta value $\zeta_{2}(a,b)$ of odd weight $k=a+b$ is a rational linear combination of 
the numbers $\zeta(k)$ and $\zeta(r)\zeta(k-r)$ where $2 \leq r \leq k/2$.} 

We illustrate the type of decomposition one obtains in Corollary \ref{DecompositionCorollary}
in the following 2-dimensional example.

\begin{example}
The double zeta function $Z_{2}\left(\F_{q}(T);s,w\right)$ has the following decomposition as
a rational linear combination of products of 1-dimensional zeta functions over $\F_{q}[T]$:
\begin{align*}
Z_{2}\left(\F_{q}(T);s,w\right)&=
\frac{q^{2}}{(q-1)^{2}}Z\left(\F_{q}[T],s+w-1\right)Z\left(\F_{q}[T],w\right)\\
&\quad -\frac{q}{(q-1)^{2}}Z\left(\F_{q}[T],s+w\right)Z\left(\F_{q}[T],w+1\right)\\
&\quad -\frac{q}{(q-1)^{2}}Z\left(\F_{q}[T],s+w\right)Z\left(\F_{q}[T],w\right)\\
&\quad +\frac{1}{(q-1)^{2}}Z\left(\F_{q}[T],s+w+1\right)Z\left(\F_{q}[T],w+1\right). 
\end{align*}
\end{example}
\vspace{0.05in}

5. If $K$ is a global function field of genus $g \geq 1$, $Z_{d}\left(K;s_{1},\ldots, s_{d}\right)$ 
no longer decomposes as a rational linear combination of 1-dimensional zeta functions. 
This indicates a more complicated arithmetic structure.
The following result establishes the analytic properties of (\ref{MDSFunctionFields})
in the case $K$ has genus $g \geq 1$. 

\begin{theorem}\label{MainTheorem3} 
Let $K$ be a global function field of genus $g \geq 1$. 

(1) $Z_{d}\left(K;s_{1},\ldots, s_{d}\right)$ is absolutely 
convergent and analytic in the region (\ref{region}). 

(2) $Z_{d}\left(K;s_{1},\ldots, s_{d}\right)$ has a meromorphic continuation 
to all $(s_{1},\ldots,s_{d})$ in $\C^{d}$ and is a rational function in each of $q^{-s_{1}}, \ldots, q^{-s_{d}}$
with a specified denominator.

(3) $Z_{d}\left(K;s_{1},\ldots, s_{d}\right)$ has possible simple poles on 
the linear subvarieties 
\begin{align*}
s_{k}+\cdots +s_{d}=0,1, \ldots, d-k+1, \quad k=1, \ldots, d.
\end{align*}
\end{theorem}

As indicated in Theorem \ref{MainTheorem3}, part $(2)$, the denominator of the rational function
$Z_{d}\left(K;s_{1},\ldots, s_{d}\right)$ can always be specified. The following result illustrates this
in 2-dimensions. 

\begin{corollary}\label{GeneralCaseCorollary} 
Let $K$ be a global function field of genus $g \geq 1$. Let $u=q^{-(s+w)}$ and $v=q^{-w}$.
Then
\begin{align*}
Z_{2}\left(K;s,w\right)=\frac{P(u,v)}{Q(u,v)},
\end{align*}
where 
$$P(u,v) \in \Q[u,v]$$ 
is a polynomial of degree $ \leq 2g+1$ in $u$ and degree $ \leq (1+2+\cdots +2g)+2g-2$ in $v$, and  
\begin{align*}
Q(u,v)=(1-u)(1-qu)(1-q^{2}u)(1-v)(1-qv)\prod_{n=0}^{2g-2}v^{n} \in \Z[u,v].
\end{align*}
\end{corollary}

\begin{remark}\label{frobenius} 
The polynomial $P(u,v)$ can be given explicitly as follows. Let 
$b_{n}$ be the number of effective divisors of $K$ of degree $n$. Let 
\begin{align*}
P_{1}(u,v)&=Q(u,v)\sum_{n=0}^{2g-2}\sum_{m=0}^{2g-2-n}b_{n}b_{m+n}u^{n}v^{m},\\
P_{2}(u,v)&=\frac{h_{K}}{q-1}(1-u)(1-qu)(1-q^{2}u)
\left[q^{g}(1-v)-(1-qv)\right]v^{2g-1}\sum_{k=0}^{2g-2}b_{k}u^{k}\prod_{\substack{n=0 \\ n \neq k}}^{2g-2}v^{n},\\
P_{3}(u,v)&=\frac{h_{K}^{2}}{(q-1)^{2}}\prod_{n=0}^{2g-2}v^{n} \quad \times \\
& \quad \left[(1-qu)(1-v)(1-u)\left(q^{2}\right)^{g}-\left(1-q^{2}u\right)(1-v)(1-u)q^{g}\right.\\
& \left. \qquad
-(1-qv)\left(1-q^{2}u\right)(1-u)q^{g}+(1-qv)\left(1-q^{2}u\right)(1-qu)\right]u^{2g-1}.
\end{align*}
Then
\begin{align*}
P(u,v)=P_{1}(u,v)+P_{2}(u,v)+P_{3}(u,v).
\end{align*}
\end{remark}

\begin{remark}\label{decomposition}
The rational function $P(u,v)/Q(u,v)$ depends in a complicated way
on the function field $K$. For example, if $K$ has genus $g=1$, then
\begin{align*}
Z_{2}(K;s,w)&=Z(K,w)+\frac{h_{K}^{2}}{(q-1)^{2}}\frac{\prod_{n=0}^{2g-2}v^{n}}{Q(u,v)} \quad \times \\
& \quad \left[(1-qu)(1-v)(1-u)q^{2}-\left(1-q^{2}u\right)(1-v)(1-u)q^{g}\right.\\
& \left. \qquad  
-(1-qv)\left(1-q^{2}u\right)(1-u)q^{g}+(1-qv)\left(1-q^{2}u\right)(1-qu)\right]u^{2g-1}. 
\end{align*}
\end{remark}
\vspace{0.05in}

6. In conclusion, we want to emphasize that $Z_{d}(K;s_{1}, \ldots, s_{d})$ provides ample
opportunity for further research. Three potentially interesting questions are the following.
First, do the 
special values $Z_{d}(K;a_{1}, \ldots a_{d})$, $a_{k} \in \Z_{ \geq 1}$, $a_{d} > 1$, 
have a Chen-type integral representation analogous to $\zeta_{d}(a_{1},\ldots,a_{d})$, and if so, does this 
lead to a cohomological interpretation of these special values?
Second, in 1-dimension the polynomial appearing in the numerator of the rational
function $Z(K,s)$ is the characteristic polynomial of the action of the Frobenius automorphism 
on the Tate module (see \cite{R}, pg. 275). Is there a similar interpretation of the polynomial $P(u,v)$? 
Third, for $K$ of genus $g \geq 1$, can  
anything in general can be said about the dependence of the rational function 
$Z_{d}(K;s_{1},\ldots, s_{d})$ on the function field $K$?

This paper is organized as follows. In section \ref{MainTheorem1Proof} we prove Theorem \ref{MainTheorem1}.
In section \ref{MainTheorem2Proof} we prove Theorem \ref{MainTheorem2} and 
Corollary \ref{DecompositionCorollary}. Finally, in section 
\ref{MainTheorem3Proof} we prove Theorem \ref{MainTheorem3} and Corollary \ref{GeneralCaseCorollary}.

I would like to thank Sol Friedberg for encouraging me to pursue this work. 
The author was supported by a Postdoctoral Fellowship at the Max-Planck-Institut f\"ur Mathematik during 
part of this work.

\section{Proof of Theorem \ref{MainTheorem1}}\label{MainTheorem1Proof}

For future reference we summarize the analytic properties of $Z(\F_{q}[T],s)$.
We refer to these as properties 1-4.

\begin{enumerate}

\item $Z(\F_{q}[T],s)$ has a meromorphic continuation to all $s$ in $\C$ and is a rational
function in $q^{-s}$. In fact, $Z(\F_{q}[T],s)=1/(1-q^{1-s})$.

\item $Z(\F_{q}[T],s)$ has a simple pole at $s=1$ with residue $1/\log(q)$.

\item The function 
\begin{align*}
\xi\left(\F_{q}[T],s\right):=\frac{q^{-s}}{1-q^{-s}}Z(\F_{q}[T],s)
\end{align*}
satisfies the functional equation
\begin{align*}
\xi\left(\F_{q}[T],1-s\right)=\xi\left(\F_{q}[T],s\right).
\end{align*}

\item $Z(\F_{q}[T],s)$ has the Euler product
\begin{align*}
Z(\F_{q}[T],s)=
\prod_{\substack{P \in \F_{q}[T]\\ P~\textrm{monic}\\ P~\textrm{irreducible}}}
\left(1-\frac{1}{\abs{P}^{s}}\right)^{-1}.
\end{align*}
\end{enumerate}

The following fact will be needed in the proof of Theorem \ref{MainTheorem1}.

\begin{lemma}\label{Monic}
The number of monic polynomials of degree $n$ in $\F_{q}[T]$ is $q^{n}$.
\end{lemma}

\textit{Proof of Theorem \ref{MainTheorem1}}. Given a set $X$, let $\abs{X}$ 
denote the number of elements in $X$. Define the nonnegative integers
\begin{align*}
a_{n_{1},\ldots, n_{d}}&=
\abs{\left\{(f_{1},\ldots,f_{d})\in \F_{q}[T]^{d}:
~f_{k}~\textrm{monic},~\textrm{deg}(f_{k})=n_{k},~k=1,\ldots,d \right\}},
\end{align*}
and 
\begin{align*}
b_{n_{k}}&=\abs{\left\{f_{k} \in \F_{q}[T]:~f_{k}~\textrm{monic},~\textrm{deg}(f_{k})=n_{k}\right\}}.
\end{align*}
Then
\begin{align*}
a_{n_{1},\ldots, n_{d}}=\prod_{k=1}^{n}b_{n_{k}},
\end{align*}
so that formally, 
\begin{align*}
Z_{d}\left(\F_{q}[T];s_{1}, \ldots, s_{d}\right)&=
\sum_{0 \leq n_{1} \leq \cdots \leq n_{d}}a_{n_{1},\ldots, n_{d}}
\prod_{k=1}^{d}\left(q^{n_{k}}\right)^{-s_{k}}\\
&=\sum_{0 \leq n_{1} \leq \cdots \leq n_{d}}
\prod_{k=1}^{d}b_{n_{k}}\left(q^{n_{k}}\right)^{-s_{k}}.
\end{align*}
The last sum can be expressed as  
\begin{align*}
\sum_{0 \leq n_{1} \leq \cdots \leq n_{d}}
\prod_{k=1}^{d}b_{n_{k}}\left(q^{n_{k}}\right)^{-s_{k}}=
\sum_{n_{1}=0}^{\infty} \cdots \sum_{n_{d}=0}^{\infty}
\prod_{k=1}^{d}b_{n_{1}+\cdots +n_{k}}\left(q^{-s_{k}}\right)^{n_{1}+\cdots + n_{k}}.
\end{align*}
Therefore,
\begin{align}\label{Decomposition}
Z_{d}\left(\F_{q}[T];s_{1},\ldots, s_{d}\right)=
\sum_{n_{1}=0}^{\infty} \cdots \sum_{n_{d}=0}^{\infty}
\prod_{k=1}^{d}b_{n_{1}+\cdots +n_{k}}\left(q^{-s_{k}}\right)^{n_{1}+\cdots + n_{k}}.
\end{align}
By Lemma \ref{Monic}, 
\begin{align*}
b_{n_{1}+\cdots +n_{k}}=q^{n_{1}+\cdots +n_{k}},
\end{align*}
so that
\begin{align*}
\prod_{k=1}^{d}b_{n_{1}+\cdots +n_{k}}\left(q^{-s_{k}}\right)^{n_{1}+\cdots + n_{k}}&=
\prod_{k=1}^{d}q^{n_{1}+\cdots +n_{k}}\left(q^{-s_{k}}\right)^{n_{1}+\cdots + n_{k}}\\
&=\prod_{k=1}^{d}\left(q^{d-k+1}\right)^{n_{k}}\left(q^{-(s_{k}+\cdots +s_{d})}\right)^{n_{k}}\\
&=\prod_{k=1}^{d}\left(q^{d-k+1-(s_{k}+\cdots +s_{d})}\right)^{n_{k}}.
\end{align*}
Thus, if $$\textrm{Re}(s_{k}+\cdots +s_{d})>d-k+1,\quad k=1,\ldots, d,$$
substituting in (\ref{Decomposition}) and summing geometric series yields
\begin{align}\label{KeyIdentity}
Z_{d}\left(\F_{q}[T];s_{1},\ldots, s_{d}\right)&=
\prod_{k=1}^{d}\sum_{n_{k}=0}^{\infty}\left(q^{d-k+1-(s_{k}+\cdots +s_{d})}\right)^{n_{k}}\notag\\
&=\prod_{k=1}^{d}\left(1-q^{d-k+1-(s_{k}+\cdots +s_{d})}\right)^{-1}.
\end{align}
This proves $(1)$.

It follows from (\ref{KeyIdentity}) that
$Z_{d}\left(\F_{q}[T];s_{1},\ldots, s_{d}\right)$ has a 
meromorphic continuation to all $(s_{1},\ldots, s_{d})$ in $\C^{d}$ and is a rational function 
in $q^{-s_{1}}, \ldots, q^{-s_{d}}$ with 
\begin{align*}
\prod_{k=1}^{d}\left(1-q^{d-k+1-(s_{k}+\cdots +s_{d})}\right)
Z_{d}\left(\F_{q}[T];s_{1},\ldots, s_{d}\right)=1.
\end{align*}  
This proves $(2)$.

It also follows from (\ref{KeyIdentity}) that 
$Z_{d}\left(\F_{q}[T];s_{1},\ldots, s_{d}\right)$ has simple poles
on the linear subvarieties
\begin{align*}
s_{k}+\cdots +s_{d}=d-k+1, \quad k=1,\ldots, d.
\end{align*}
This proves $(3)$.

Write
\begin{align*}
d-k+1-(s_{k}+\cdots +s_{d})=1-(s_{k}+\cdots +s_{d}-(d-k)).
\end{align*}
Then (\ref{KeyIdentity}) and property 1 yield the factorization
\begin{align}\label{Factorization1}
Z_{d}\left(\F_{q}[T];s_{1},\ldots, s_{d}\right)=
\prod_{k=1}^{d}Z\left(\F_{q}[T],s_{k}+\cdots +s_{d}-(d-k)\right).
\end{align}
From the definition of $\xi(\F_{q}[T],s)$ and (\ref{Factorization1}) we obtain 
\begin{align}\label{Factorization2}
\xi_{d}\left(\F_{q}[T];s_{1},\ldots, s_{d}\right)
&=\prod_{k=1}^{d}\xi\left(\F_{q}[T],s_{k}+\cdots +s_{d}-(d-k)\right).
\end{align}
From property 3 we obtain the functional equation
\begin{align*}
\xi\left(\F_{q}[T],1-s_{1}+\cdots +s_{d}-(d-1)\right)=
\xi\left(\F_{q}[T],s_{1}+\cdots +s_{d}-(d-1)\right),
\end{align*}
or equivalently,
\begin{align}\label{FunctionalEquation2}
\xi\left(\F_{q}[T],d-(s_{1}+\cdots +s_{d})\right)=
\xi\left(\F_{q}[T],s_{1}+\cdots +s_{d}-(d-1)\right).
\end{align}
Let 
\begin{align*}
s_{1} \mapsto -s_{1}-2(s_{2}+\cdots +s_{d})+2d-1
\end{align*}
in 
\begin{align*}
\xi_{d}\left(\F_{q}[T];s_{1},\ldots, s_{d}\right).
\end{align*}
It follows that
\begin{align*}
&\xi_{d}\left(\F_{q}[T];-s_{1}-2(s_{2}+\cdots +s_{d})+2d-1,\ldots, s_{d}\right)\\
&=\xi\left(\F_{q}[T],\left(-s_{1}-2(s_{2}+\cdots +s_{d})+2d-1\right)+
s_{2}+\cdots +s_{d}-(d-1)\right)\\
& \quad \times \prod_{k=2}^{d}\xi\left(\F_{q}[T],s_{k}+\cdots +s_{d}-(d-k)\right)\\
&=\xi\left(\F_{q}[T],d-(s_{1}+\cdots +s_{d})\right)
\prod_{k=2}^{d}\xi\left(\F_{q}[T],s_{k}+\cdots +s_{d}-(d-k)\right)\\
&=\xi\left(\F_{q}[T],s_{1}+\cdots +s_{d}-(d-1)\right)
\prod_{k=2}^{d}\xi\left(\F_{q}[T],s_{k}+\cdots +s_{d}-(d-k)\right) \quad 
(\textrm{from (\ref{FunctionalEquation2}))}\\
&=\prod_{k=1}^{d}\xi\left(\F_{q}[T],s_{k}+\cdots +s_{d}-(d-k)\right)\\
&=\xi_{d}\left(\F_{q}[T];s_{1},\ldots, s_{d}\right) \quad 
(\textrm{from (\ref{Factorization2}))}.
\end{align*}
This proves $(4)$.

Finally, from the factorization (\ref{Factorization1}) and property 4 
we obtain the Euler product
\begin{align*}
Z_{d}\left(\F_{q}[T];s_{1},\ldots, s_{d}\right)=
\prod_{\substack{P \in \F_{q}[T]\\ P~\textrm{monic}\\P~\textrm{irreducible}}}
\prod_{k=1}^{d}\left(1-\frac{1}{\abs{P}^{s_{k}+\cdots +s_{d}-(d-k)}}\right)^{-1}.
\end{align*}
This proves $(5)$.
\qed

\section{Proofs of Theorem \ref{MainTheorem2} and Corollary \ref{DecompositionCorollary}}\label{MainTheorem2Proof}
It can be shown that for all integers $n \geq 0$, the number $b_{n}$ of effective divisors of degree $n$ is finite 
(see \cite{R}, Lemma 5.5). It is known that if $n$ is sufficiently large compared to the genus of $K$, this number
takes an explicit form. 

\begin{proposition}\label{EffectiveNumber}
For all nonnegative integers $n > 2g-2$, 
\begin{align*}
b_{n}=h_{K}\frac{q^{n-g+1}-1}{q-1}.
\end{align*}
\end{proposition}

We include a proof for convenience.

\begin{proof} Let $A$ be a divisor and $\overline{A}$ be its divisor class. We will need
the following two facts. 

\begin{enumerate}
\item The number of effective divisors in $\overline{A}$ is 
\begin{align*}
\frac{q^{l(A)}-1}{q-1}.
\end{align*}
\item If $\deg(A) > 2g-2$, then $l(A)=\deg(A)-g+1$.
\end{enumerate}

Fact 1 is \cite{R}, Lemma 5.7. 

To prove Fact 2, first observe that by Theorem \ref{RiemannRoch}, Fact 2 is equivalent to $l(C-A)=0$. 
Now, by Theorem \ref{RiemannRoch} with $A=C$, 
$\deg(C)=l(C)+g-2$ (here we used that $l(0)=1$), and by Theorem \ref{RiemannRoch} with $A=0$, 
$l(C)=g$. Thus, $\deg(C)=2g-2.$ Since we have assumed that $\deg(A) > 2g-2$, it follows that 
\begin{align*}
\deg(C-A)=\deg(C)-\deg(A)<2g-2+2-2g=0.
\end{align*}
Finally, by \cite{R}, Lemma 5.3, $\deg(C-A) < 0$ implies that $l(C-A)=0$. This proves Fact 2. 

We can now prove the proposition. From the exact sequence
\begin{align*}
0 \rightarrow Cl_{K}^{0} \rightarrow Cl_{K} \rightarrow \Z \rightarrow 0
\end{align*}
we conclude that for each nonnegative integer $n$, there are $h_{K}=\abs{Cl_{K}^{0}}$ divisor classes of degree $n$. 
List these as $\left\{\overline{A}_{1}, \ldots, \overline{A}_{h_{K}}\right\}$. By Fact 1, 
the number of effective divisors in $\overline{A_{i}}$ is 
\begin{align*}
\frac{q^{l(A_{i})}-1}{q-1}.
\end{align*}
Therefore,
\begin{align*}
b_{n}=\sum_{i=1}^{h_{K}}\frac{q^{l(A_{i})}-1}{q-1}.
\end{align*}
Assume that $n=\deg(A_{i}) > 2g-2$. By Fact 2, $l(A_{i})=n-g+1$. We conclude that
\begin{align*}
b_{n}=h_{K}\frac{q^{n-g+1}-1}{q-1}.
\end{align*}
\end{proof}

\textit{Proof of Theorem \ref{MainTheorem2}.}
Define the nonnegative integers
\begin{align*}
a_{n_{1},\ldots, n_{d}}&=
\abs{\left\{(D_{1},\ldots,D_{d})\in D_{K}^{+} \times \cdots \times D_{K}^{+}:
~\textrm{deg}(D_{k})=n_{k},~k=1,\ldots,d \right\}},
\end{align*}
and 
\begin{align*}
b_{n_{k}}&=\abs{\left\{D_{k} \in D_{K}^{+}:~\textrm{deg}(D_{k})=n_{k}\right\}}.
\end{align*}
Then
\begin{align*}
a_{n_{1},\ldots, n_{d}}=\prod_{k=1}^{n}b_{n_{k}},
\end{align*}
so that formally,
\begin{align}\label{Decomposition2}
Z_{d}\left(\F_{q}(T);s_{1},\ldots, s_{d}\right)=
\sum_{n_{1}=0}^{\infty} \cdots \sum_{n_{d}=0}^{\infty}
\prod_{k=1}^{d}b_{n_{1}+\cdots +n_{k}}\left(q^{-s_{k}}\right)^{n_{1}+\cdots + n_{k}}.
\end{align}
By Proposition \ref{EffectiveNumber}, 
\begin{align*}
\abs{b_{n_{1}+\cdots +n_{k}}} \leq C_{k} q^{n_{1}+\cdots +n_{k}}
\end{align*}
for some constants $C_{k} > 0$, $k=1, \ldots, d$. This yields the estimate
\begin{align*}
\prod_{k=1}^{d}\abs{b_{n_{1}+\cdots +n_{k}}\left(q^{-s_{k}}\right)^{n_{1}+\cdots + n_{k}}} &\leq
\prod_{k=1}^{d}C_{k}\abs{q^{n_{1}+\cdots n_{k}} \left(q^{-s_{k}}\right)^{n_{1}+\cdots + n_{k}}}\\
&=\prod_{k=1}^{d}C_{k}\left(q^{d-k+1-\textrm{Re}(s_{k}+\cdots +s_{d})}\right)^{n_{k}}.
\end{align*}
Thus, if $$\textrm{Re}(s_{k}+\cdots +s_{d})>d-k+1, \quad k=1,\ldots, d,$$
substituting in (\ref{Decomposition2}) and summing geometric series yields
\begin{align*}
\sum_{n_{1}=0}^{\infty}\cdots \sum_{n_{d}=0}^{\infty}
\prod_{k=1}^{d}\abs{b_{n_{1}+\cdots +n_{k}}\left(q^{-s_{k}}\right)^{n_{1}+\cdots + n_{k}}}
& \leq \prod_{k=1}^{d}C_{k}\sum_{n_{k}=0}^{\infty}q^{d-k+1-\textrm{Re}(s_{k}+\cdots +s_{d})}\\
&=\prod_{k=1}^{d}C_{k}\left(1-q^{d-k+1-\textrm{Re}(s_{k}+\cdots +s_{d})}\right)^{-1}.
\end{align*}
This proves $(1)$ for any global function field.

Suppose first that $d=2$. Because $\F_{q}(T)$ has genus $g=0$ and class number $h_{\F_{q}(T)}=1$, 
Proposition \ref{EffectiveNumber} implies that 
\begin{align*}
b_{n}=\frac{q^{n+1}-1}{q-1}
\end{align*}
for all integers $n \geq 0$. Then
\begin{align*}
&b_{n}b_{n+m}\left(q^{-s}\right)^{n}\left(q^{-w}\right)^{n+m}\notag\\
&=\frac{1}{(q-1)^{2}}\left(q^{n+1}-1\right)\left(q^{n+m+1}-1\right)
\left(q^{-s}\right)^{n}\left(q^{-w}\right)^{n+m}\notag\\
&=\frac{1}{(q-1)^{2}}\left[q^{n+1}q^{n+m+1}-q^{n+1}-q^{n+m+1}+1\right]
\left(q^{-(s+w)}\right)^{n}\left(q^{-w}\right)^{m}
\end{align*}
for all integers $n, m \geq 0$. 
Substitute in (\ref{Decomposition2}) and expand to obtain
\begin{align*}
&(q-1)^{2}Z_{2}\left(\F_{q}(T);s,w\right)\\
&=q^{2}\sum_{n=0}^{\infty}\left(q^{2-(s+w)}\right)^{n}
\sum_{m=0}^{\infty}\left(q^{1-w}\right)^{m} 
-q\sum_{n=0}^{\infty}\left(q^{1-(s+w)}\right)^{n}
\sum_{m=0}^{\infty}\left(q^{-w}\right)^{m}\\
& \quad -q\sum_{n=0}^{\infty}\left(q^{1-(s+w)}\right)^{n}
\sum_{m=0}^{\infty}\left(q^{1-w}\right)^{m}
+\sum_{n=0}^{\infty}\left(q^{-(s+w)}\right)^{n}
\sum_{m=0}^{\infty}\left(q^{-w}\right)^{m}. 
\end{align*}
Summing geometric series in the preceding expression yields
\begin{align}\label{RationalExpression}
&Z_{2}\left(\F_{q}(T);s,w\right)=\\
&\quad \frac{q^{2}}{(q-1)^{2}}\frac{1}{1-q^{2-(s+w)}}\frac{1}{1-q^{1-w}}
-\frac{q}{(q-1)^{2}}\frac{1}{1-q^{1-(s+w)}}\frac{1}{1-q^{-w}}\notag\\
&\quad -\frac{q}{(q-1)^{2}}\frac{1}{1-q^{1-(s+w)}}\frac{1}{1-q^{1-w}}
+\frac{1}{(q-1)^{2}}\frac{1}{1-q^{-(s+w)}}\frac{1}{1-q^{-w}}\notag.
\end{align}

It follows from (\ref{RationalExpression}) that $Z_{2}\left(\F_{q}(T);s,w\right)$ has a 
meromorphic continuation to all $(s,w)$ in $\C^{2}$ and is a rational function in $q^{-s}$ and 
$q^{-w}$.

Let
\begin{align*}
&Q\left(q^{-s},q^{-w}\right)=\\
&\quad (q-1)^{2}\left(1-q^{-(s+w)}\right)\left(1-q^{1-(s+w)}\right)\left(2-q^{-(s+w)}\right)
\left(1-q^{-w}\right)\left(1-q^{1-w)}\right).
\end{align*}
Multiply both sides of (\ref{RationalExpression}) by $Q\left(q^{-s},q^{-w}\right)$ to obtain
\begin{align*}
&Q\left(q^{-s},q^{-w}\right)Z_{2}\left(\F_{q}(T);s,w\right)=\\
&\quad q^{2}\left(1-q^{-(s+w)}\right)\left(1-q^{1-(s+w)}\right)\left(1-q^{-w}\right)\notag\\
& \quad -q\left(1-q^{-(s+w)}\right)\left(1-q^{2-(s+w)}\right)\left(1-q^{1-w}\right)\notag\\
&\quad -q\left(1-q^{-(s+w)}\right)\left(1-q^{2-(s+w)}\right)\left(1-q^{-w}\right)\notag\\
& \quad +\left(1-q^{1-(s+w)}\right)\left(1-q^{2-(s+w)}\right)\left(1-q^{1-w}\right).\notag
\end{align*}
Thus,  
$$Q\left(q^{-s},q^{-w}\right)Z_{2}\left(\F_{q}(T);s,w\right)$$
is polynomial of degree $\leq 3$ in each of $q^{-s}$ and $q^{-w}$. This proves $(2)$.

It also follows from (\ref{RationalExpression}) that $Z_{2}\left(\F_{q}(T);s,w\right)$ has possible
simple poles on the linear subvarieties $s+w=0,1,2$ and $w=0,1$. This proves $(3)$. 

The proof for $d\geq 3$ is analogous. Substitute the product
\begin{align*}
\prod_{k=1}^{d}b_{n_{1}+\cdots +n_{k}}=\frac{1}{(q-1)^{d}}
\prod_{k=1}^{d}\left(q^{n_{1}+\cdots +n_{k}+1}-1\right)
\end{align*} 
into (\ref{Decomposition2}),  
expand, and sum geometric series to obtain 
\begin{align}\label{RationalExpressionSecond}
Z_{d}\left(\F_{q}(T);s_{1}, \ldots, s_{d}\right)=
\sum_{i=1}^{2^{d}}R_{i}\left(q^{-s_{1}}, \ldots, q^{-s_{d}}\right),
\end{align}
where $R_{i}\left(q^{-s_{1}}, \ldots, q^{-s_{d}}\right)$ is a rational function which is a product
of one function from each of the sets
\begin{align*}
&\left\{(q-1)^{-3}, \ldots, (q-1)^{-d}\right\},\\
&\left\{q, \ldots, q^{d}\right\},\\
&\left\{\left(1-q^{-s_{d}}\right)^{-1}, 
\left(1-q^{1-s_{d}}\right)^{-1}\right\},
\end{align*}
and $d-1$ functions from the set
\begin{align*}
\left\{\left(1-q^{-(s_{k}+\cdots + s_{d})}\right)^{-1},~\left(1-q^{1-(s_{k}+\cdots + s_{d})}\right)^{-1},~ 
\left(1-q^{2-(s_{k}+\cdots + s_{d})}\right)^{-1}:~k=1, \ldots, d-1.\right\}
\end{align*}

It follows from (\ref{RationalExpressionSecond}) that 
$Z_{d}\left(\F_{q}(T);s_{1}, \ldots, s_{d}\right)$ has a meromorphic continuation to
all $(s_{1},\ldots, s_{d})$ in $\C^{d}$ and is a rational function in $q^{-s_{1}}, \ldots,
q^{-s_{d}}$. 

Let 
\begin{align*}
&Q\left(q^{-s_{1}},\ldots, q^{-s_{d}}\right)
=(q-1)^{d}\left(1-q^{-s_{d}}\right)\left(1-q^{1-s_{d}}\right)\\
& \times \prod_{k=1}^{d-1}\left(1-q^{-(s_{k}+\cdots +s_{d})}\right)
\left(1-q^{1-(s_{k}+\cdots +s_{d})}\right)
\left(1-q^{2-(s_{k}+\cdots +s_{d})}\right).
\end{align*}
Then it is not difficult to show that 
\begin{align*}
Q\left(q^{-s_{1}},\ldots, q^{-s_{d}}\right)R_{i}\left(q^{-s_{1}}, \ldots, q^{-s_{d}}\right)
\end{align*}
is a polynomial of degree $ \leq 2d-1$ in each of $q^{-s_{1}}, \ldots, q^{-s_{d}}$.

Finally, it follows from the explicit form of the functions 
$R_{i}\left(q^{-s_{1}}, \ldots, q^{-s_{d}}\right)$ that 
\begin{align*}
\sum_{i=1}^{2^{d}}R_{i}\left(q^{-s_{1}}, \ldots, q^{-s_{d}}\right)
\end{align*}
has possible simple poles on the linear subvarieties
\begin{align*}
s_{k}+\cdots +s_{d}=0,1, \ldots, d-k+1, \quad k=1, \ldots, d.
\end{align*}
\qed
\vspace{0.10in}

\textit{Proof of Corollary \ref{DecompositionCorollary}.} From property 1 we see that each
function $R_{i}\left(q^{-s_{1}}, \ldots, q^{-s_{d}}\right)$ is the product of a rational 
number and products of zeta functions of the type 
\begin{align*}
Z\left(\F_{q}[T],s_{k}+\cdots +s_{d}-1\right),~Z\left(\F_{q}[T],s_{k}+\cdots +s_{d}\right),~
Z\left(\F_{q}[T],s_{k}+\cdots +s_{d}+1\right),
\end{align*}
where $k=1, \ldots, d$. The corollary now follows from (\ref{RationalExpressionSecond}).
\qed

\section{Proof of Theorem \ref{MainTheorem3} and Corollary \ref{GeneralCaseCorollary}}
\label{MainTheorem3Proof}

\textit{Proof of Theorem \ref{MainTheorem3}.}
Part $(1)$ was established in the proof of Theorem \ref{MainTheorem2}. 

We will prove parts $(2)$ and $(3)$ for $d=2$, the proofs for $d \geq 3$ being a more complicated
elaboration on the same idea.

Write 
\begin{align*}
b_{n}b_{n+m}\left(q^{-s}\right)^{n}\left(q^{-w}\right)^{n+m}=
b_{n}b_{n+m}\left(q^{-(s+w)}\right)^{n}\left(q^{-w}\right)^{m}
\end{align*}
and substitute in (\ref{Decomposition2}) to obtain 
\begin{align}\label{Decomposition3}
Z_{2}\left(K;s,w\right)=\sum_{n=0}^{\infty}b_{n}\left(q^{-(s+w)}\right)^{n}
\sum_{m=0}^{\infty}b_{n+m}\left(q^{-w}\right)^{m}.
\end{align}

Let $u=q^{-(s+w)}$ and $v=q^{-w}$. Decompose the sum (\ref{Decomposition3}) as follows:
\begin{align*}
&\sum_{n=0}^{\infty}b_{n}u^{n}\sum_{m=0}^{\infty}b_{n+m}v^{m}\\
&=\sum_{n=0}^{2g-2}b_{n}u^{n}\sum_{m=0}^{\infty}b_{n+m}v^{m}+
\sum_{n=2g-1}^{\infty}b_{n}u^{n}\sum_{m=0}^{\infty}b_{n+m}v^{m}\\
&=\sum_{n=0}^{2g-2}b_{n}u^{n}\sum_{m=0}^{2g-2-n}b_{m+n}v^{m}+
\sum_{n=0}^{2g-2}b_{n}u^{n}\sum_{m=2g-1-n}^{\infty}b_{m+n}v^{m}+
\sum_{n=2g-1}^{\infty}b_{n}u^{n}\sum_{m=0}^{\infty}b_{n+m}v^{m}\\
&=A(u,v)+B(u,v)+C(u,v).
\end{align*}

It is immediate that 
\begin{align*}
A(u,v)=\sum_{n=0}^{2g-2}\sum_{m=0}^{2g-2-n}b_{n}b_{m+n}u^{n}v^{m}  
\end{align*}
is analytic for all $(s,w)$ in $\C^{2}$ and is a rational function in $q^{-s}$ and $q^{-w}$. 

To analyze $B(u,v)$, first observe that if  $m \geq 2g-1-n$, then 
$m+n>2g-2$. Applying Proposition
\ref{EffectiveNumber} and summing geometric series yields
\begin{align*}
\sum_{m=2g-1-n}^{\infty}b_{m+n}v^{m}&=
\frac{h_{K}}{q-1}\sum_{m=2g-1-n}^{\infty}\left(q^{m+n-g+1}-1\right)v^{m}\\
&=\frac{h_{K}}{q-1}\left[q^{n-g+1}\sum_{m=2g-1-n}^{\infty}(qv)^{m}-
\sum_{m=2g-1-n}^{\infty}v^{m}\right]\\
&=\frac{h_{K}}{q-1}\left[q^{n-g+1}\frac{(qv)^{2g-1-n}}{1-qv}-\frac{v^{2g-1-n}}{1-v}\right]\\
&=\frac{h_{K}}{q-1}\left[\frac{q^{g}}{1-qv}-\frac{1}{1-v}\right]v^{2g-1-n}.
\end{align*}
Substitute this computation into $B(u,v)$ to obtain
\begin{align*}
B(u,v)=\frac{h_{K}}{q-1}\left[\frac{q^{g}}{1-qv}-\frac{1}{1-v}\right]v^{2g-1}
\sum_{n=0}^{2g-2}b_{n}\left(uv^{-1}\right)^{n}.
\end{align*}

This expression shows that $B(u,v)$ has a meromorphic continuation to all $(s,w)$ in 
$\C^{2}$ and is a rational function in $q^{-s}$ and $q^{-w}$. 
Further, this expression shows that $B(u,v)$ has simple poles at  
$v=1$ and $v=q^{-1}$, which correspond to the linear subvarieties $w=0$ and $w=1$, respectively.

To analyze $C(u,v)$, first observe that if $n \geq 2g-1$, then $m+n > 2g-2$ for all $m \geq 0$.
Applying Proposition \ref{EffectiveNumber} and summing geometric series yields
\begin{align*}
\sum_{m=0}^{\infty}b_{m+n}v^{m}&=\frac{h_{K}}{q-1}\sum_{m=0}^{\infty}\left(q^{m+n-g+1}-1\right)v^{m}\\
&=\frac{h_{K}}{q-1}\left[q^{n-g+1}\sum_{m=0}^{\infty}(qv)^{m}-\sum_{m=0}^{\infty}v^{m}\right]\\
&=\frac{h_{K}}{q-1}\left[q^{n-g+1}\frac{1}{1-qv}-\frac{1}{1-v}\right].
\end{align*} 
Substitute this computation into $C(u,v)$, apply Proposition \ref{EffectiveNumber}, 
and sum geometric series to obtain 
\begin{align*}
&C(u,v)\\
&=\frac{h_{K}}{q-1}\left[\frac{q^{-g+1}}{1-qv}\sum_{n=2g-1}^{\infty}b_{n}(qu)^{n}-
\frac{1}{1-v}\sum_{n=2g-1}^{\infty}b_{n}u^{n}\right]\\
&=\frac{h_{K}}{q-1} \quad \times \\
&\quad \left[\frac{q^{-g+1}}{1-qv}\sum_{n=2g-1}^{\infty}
\left[\frac{h_{K}}{q-1}\left(q^{n-g+1}-1\right)\right](qu)^{n}
-\frac{1}{1-v}\sum_{n=2g-1}^{\infty}
\left[\frac{h_{K}}{q-1}\left(q^{n-g+1}-1\right)\right]u^{n}\right]\\
&=\left(\frac{h_{K}}{q-1}\right)^{2}\quad \times \\
&\quad \left[\frac{q^{-g+1}}{1-qv}
\left[q^{-g+1}\sum_{n=2g-1}^{\infty}(q^{2}u)^{n}-\sum_{n=2g-1}^{\infty}(qu)^{n}\right]
-\frac{1}{1-v}\left[q^{-g+1}\sum_{n=2g-1}^{\infty}(qu)^{n}
-\sum_{n=2g-1}^{\infty}u^{n}\right]\right]\\
&=\left(\frac{h_{K}}{q-1}\right)^{2} \quad \times \\
&\quad \left[\frac{q^{-g+1}}{1-qv}
\left[q^{-g+1}\frac{(q^{2}u)^{2g-1}}{1-q^{2}u}-\frac{(qu)^{2g-1}}{1-qu}\right]
-\frac{1}{1-v}\left[q^{-g+1}\frac{(qu)^{2g-1}}{1-qu}-\frac{u^{2g-1}}{1-u}\right]\right]\\
&=\left(\frac{h_{K}}{q-1}\right)^{2} \quad \times \\
&\quad \left[\frac{(q^{2})^{g}}{(1-qv)(1-q^{2}u)}
-\frac{q^{g}}{(1-qv)(1-qu)}-\frac{q^{g}}{(1-v)(1-qu)}+\frac{1}{(1-v)(1-u)}\right]u^{2g-1}.
\end{align*}

This expression shows that $C(u,v)$ has a meromorphic continuation to all $(s,w)$ in $\C^{2}$ 
and is a rational function in $q^{-s}$ and $q^{-w}$. 
Further, this expression shows that $C(u,v)$ has possible simple poles at
$v=1$, $v=q^{-1}$, $u=1$, $u=q^{-1}$, and $u=q^{-2}$, which correspond to the linear
subvarieties $w=0$, $w=1$, $s+w=0$, $s+w=1$, and $s+w=2$, respectively.  \qed
\vspace{0.10in}

\textit{Proof of Corollary \ref{GeneralCaseCorollary}.}
First observe that the polynomial 
\begin{align*}
Q(u,v)=(1-qv)\left(1-q^{2}u\right)(1-qu)(1-v)(1-u)\prod_{n=0}^{2g-2}v^{n}
\end{align*}
has degree 3 in $u$ and degree $1+2+\cdots +2g$ in $v$. 

Because 
\begin{align*}
A(u,v)=\sum_{n=0}^{2g-2}\sum_{m=0}^{2g-2-n}b_{n}b_{m+n}u^{n}v^{m}  
\end{align*}
has degree $ \leq 2g-2$ in $u$ and degree $ \leq 2g-2$ in $v$, $Q(u,v)A(u,v)$ has degree
$\leq 2g+1$ in $u$ and degree $\leq (1+2+\cdots +2g)+2g-2$ in $v$.

Write
\begin{align*}
\sum_{n=0}^{2g-2}b_{n}\left(uv^{-1}\right)^{n}=
\frac{\sum_{k=0}^{2g-2}b_{k}u^{k}\prod_{\substack{n=0\\n \neq k}}^{2g-2}v^{n}}
{\prod_{n=0}^{2g-2}v^{n}},
\end{align*}
so that 
\begin{align*}
B(u,v)=\frac{h_{K}}{q-1}\left[\frac{q^{g}(1-v)-(1-qv)}{(1-qv)(1-v)}\right]v^{2g-1}
\frac{\sum_{k=0}^{2g-2}b_{k}u^{k}\prod_{\substack{n=0\\n \neq k}}^{2g-2}v^{n}}
{\prod_{n=0}^{2g-2}v^{n}}.
\end{align*}
Then
\begin{align*}
&Q(u,v)B(u,v)=\\
&\frac{h_{K}}{q-1}\left(1-q^{2}u\right)(1-qu)(1-u)\left[q^{g}(1-v)-(1-qv)\right]v^{2g-1}
\sum_{k=0}^{2g-2}b_{k}u^{k}\prod_{\substack{n=0\\n \neq k}}^{2g-2}v^{n}
\end{align*}
has degree $\leq 2g+1$ in $u$ and degree $(1+2+\cdots + 2g-2)+2g$ in $v$.

Similarly, write
\begin{align*}
&C(u,v)=\frac{h_{K}^{2}}{(q-1)^{2}}\frac{1}{(1-qv)\left(1-q^{2}u\right)(1-qu)(1-v)(1-u)}\\
\quad &\times \left[(1-qu)(1-v)(1-u)q^{2}-\left(1-q^{2}u\right)(1-v)(1-u)q^{g}\right.\\
& \left. \quad \quad 
-(1-qv)\left(1-q^{2}u\right)(1-u)q^{g}+(1-qv)\left(1-q^{2}u\right)(1-qu)\right]u^{2g-1}. 
\end{align*}
Then
\begin{align*}
&Q(u,v)C(u,v)=\frac{h_{K}^{2}}{(q-1)^{2}}\prod_{n=0}^{2g-2}v^{n}\\
\quad & \times \left[(1-qu)(1-v)(1-u)\left(q^{2}\right)^{g}-\left(1-q^{2}u\right)(1-v)(1-u)q^{g}\right.\\
& \left. \quad \quad 
-(1-qv)\left(1-q^{2}u\right)(1-u)q^{g}+(1-qv)\left(1-q^{2}u\right)(1-qu)\right]u^{2g-1}
\end{align*}
has degree $2g+1$ in $u$ and degree $(1+2+\cdots +2g-2)+1$ in $v$.

The corollary follows by comparing the degrees of the polynomials on the RHS of 
the expression 
\begin{align*}
Q(u,v)Z_{2}\left(K;s,w\right)=Q(u,v)A(u,v)+Q(u,v)B(u,v)+Q(u,v)C(u,v).
\end{align*}
\qed

\end{document}